\documentclass{article}
\usepackage{latexsym}
\usepackage{amsmath,amsthm}
\usepackage{amssymb}
\usepackage{graphics}
\usepackage{color}
\newtheorem{proposition}{Proposition}
\newtheorem{remark}{Remark}
\newtheorem{theorem}{Theorem}

\def\pa{\partial}

\begin{document}
\title{Small-amplitude nonlinear waves
on a black hole background}
\author{Mihalis Dafermos\thanks{University of Cambridge,
Department of Pure Mathematics
and Mathematical Statistics, Wilberforce Road, Cambridge CB3 0WB, United 
Kingdom} \and  
Igor Rodnianski\thanks{Princeton University,
Department of Mathematics, Washington Road, Princeton NJ 08544, United 
States}}
\maketitle
\begin{abstract}
Let $G(x)$ be a $C^0$ function such that $|G(x)|\le K|x|^{p}$
for $|x|\le c$, for constants $K,c>0$.
We consider spherically symmetric
solutions of $\Box_g\phi=G(\phi)$ where $g$ is a Schwarzschild
or more generally a Reissner-Nordstr\"om metric, 
and such that $\phi$ and $\nabla \phi$ are 
compactly supported on a complete Cauchy surface.
It is proven that for $p> 4$,
such solutions do not blow up in the 
domain of outer communications, provided the initial data
are small. Moreover, $|\phi|\le C(\max\{v,1\})^{-1}$, where $v$ denotes an
Eddington-Finkelstein advanced time coordinate.
\end{abstract}
The interaction between the geometry
of black hole backgrounds and the behaviour (in the large) of
linear and non-linear waves
plays a fundamental role in general relativity. The very stability
of the simplest black hole solutions like Schwarzschild and Kerr could 
depend on subtle features of this interaction. 
Moreover, fine aspects of the behaviour of these waves, 
such as the phenomenon of decaying tails, 
have physical importance, 
both from the point of view of
far away observers, 
and for those brave enough to cross the ``event horizon''.

In~\cite{mi:mazi}, we studied the problem of
collapse of a spherically symmetric
self-gravitating scalar field, i.e.~solutions of the system
\begin{equation}
\label{sg01}
R_{\mu\nu}-\frac12g_{\mu\nu}R=2T_{\mu\nu}
\end{equation}
\begin{equation}
\label{sg02}
\Box_g\phi=0
\end{equation}
\begin{equation}
\label{sg03}
T_{\mu\nu}=\phi_{,\mu}\phi_{,\nu}-\frac12g_{\mu\nu}\phi^{,\lambda}\phi_{,\lambda},
\end{equation}
arising
 from spherically symmetric
initial data, in the case where it is known \emph{a priori}
that a black hole forms.\footnote{For the definition of the notion of
black hole, see~\cite{md:sssts}. Black holes form for generic 
large initial data in view of Christodoulou's proof of weak
cosmic censorship~\cite{chr:ins} for spherically symmetric
solutions to \eqref{sg01}--\eqref{sg03}. 
In the presence of an additional non-trivial Maxwell field, for which
the results of~\cite{mi:mazi} also apply, 
the existence of a black hole follows from~\cite{md:sssts}.}
The main result of~\cite{mi:mazi} was a proof of 
\emph{Price's law}, heuristically derived
in~\cite{rpr:ns, gpp:de1}. The mathematical content 
of this ``law'' is a set of decay rates for the scalar field in appropriately defined
null coordinates on the so-called \emph{domain of outer 
communications} or \emph{exterior} of the
black hole. In particular, we showed decay 
$|\phi|\le v^{-3+\epsilon }$ in an Eddington-Finkelstein-like 
advanced time $v$-coordinate along the event horizon, and 
decay $|r\phi|\le u^{-2}$ in retarded time $u$ along null infinity.

The problem of the self-gravitating scalar field studied
in~\cite{mi:mazi} is non-linear. In particular, the black hole
geometry depends on the solution, and even properly identifying 
the black hole exterior is a non-trivial part of the problem. 
The results of~\cite{mi:mazi} can be specialized, however,
to the easier problem of the study of spherically symmetric
solutions to the equation
\[
\Box_g\phi=0
\]
on a fixed \emph{Reissner-Nordstr\"om} 
background spacetime $(\mathcal{M},g)$. 
These spacetimes constitute a two-parameter 
family of solutions to the Einstein-Maxwell equations, with 
parameters mass $M$ and charge $e$. The one-parameter subfamily
of vacuum solutions defined by $e=0$ is known as the \emph{Schwarzschild} family.
For parameter values $0\le |e| <M$, these spacetimes contain
black and white holes, and the so-called exterior region 
can be covered by coordinates $(r,t)$ ranging in 
$(M+\sqrt{M^2-e^2},\infty)\times(-\infty, \infty)$ such that
the metric takes the form: 
\[
g=-\left(1-\frac {2M}r+\frac{e^2}{r^2}\right) dt^2 +
\left(1-\frac {2M}r+\frac{e^2}{r^2}\right)^{-1} dr^2+ 
r^2 d\sigma_{{\mathbb S}^2},
\]
where $d\sigma_{{\mathbb S}^2}$ denotes the standard metric on the unit two-sphere.
In~\cite{mi:mazi}, decay rates were obtained 
both in the (static) $t$ coordinate, as well as 
in Eddington-Finkelstein retarded and advanced time coordinates $u$ and $v$. (The latter
are indispensible for understanding the decay of the radiation tail along
null infinity and the behaviour along the event horizon.)
Previously, decay rates had been shown only in the static $t$ 
coordinate~\cite{clsy:wpgs, st:pc}
under the additional restrictive assumption that the solution is
supported away from the event horizon's sphere of bifurcation.\footnote{We also
note earlier results showing boundedness~\cite{kw:lss} and decay without a rate~\cite{ft:td}.}
In contrast to the situation in Minkowski space, there is no 
decay rate in $t$ which is uniform in $r$, 
in view of the geometry of the horizon.

In this paper, we shall explore the implications of the decay 
proven in~\cite{mi:mazi} for the behaviour of small data solutions
to the semi-linear wave equation
\begin{equation}
\label{nlw}
\Box_g\phi=G(\phi)
\end{equation}
on a fixed Reissner-Nordstr\"om background $(\mathcal{M},g)$,
where the nonlinearity is assumed to decay 
sufficiently fast at $0$:
\begin{equation}
\label{satisfies}
|G(x)|\le K|x|^{p}{\rm\ for\ }|x|\le c
\end{equation}
for some constants $K,c>0$. The study of small-data solutions to
$(\ref{nlw})$ under assumption $(\ref{satisfies})$
is often referred to as the ``John problem''. The 
point of $(\ref{satisfies})$ is that no special structure 
is assumed of the non-linearity; in particular, it is \emph{not} assumed that
 $G(x)=V'(x)$ for a nonnegative potential $V\ge 0$. Our main result (Theorem~\ref{main}) 
states that for $p>4$, small-data spherically symmetric solutions 
of~\eqref{nlw} remain regular on the closure of the exterior region
$\overline{J^-(\mathcal{I}^+)\cap J^+(\mathcal{I}^-)}$ of Schwarzschild
or Reissner Nordstr\"om.\footnote{In stating results, we 
apply here standard notations from global Lorentzian geometry. See Hawking
and Ellis~\cite{he:lssst}. The reader unfamiliar with these global concepts
should note that our result implies in particular ``global existence''
with respect to Schwarzschild or Reissner-Nordstr\"om time $t$.}
Moreover, $|\phi|\le Cv^{-1}$ in advanced Eddington-Finkelstein 
time $v$, and, for $r\ge r_0>M+\sqrt{M^2-e^2}$, $|\phi|(r,t)\le C(r_0)t^{-1}$, but
with $C\to\infty$ as $r_0\to M+\sqrt{M^2-e^2}$! Note that our 
assumptions on the initial data in Theorem~\ref{main} are the 
geometrically appropriate ones; no special conditions are imposed
at the sphere of bifurcation $\mathcal{H}^+_A\cap\mathcal{H}^+_B$ of the event horizon. 

The global study of the equation $(\ref{nlw})$
on a Schwarzschild background was initiated in~\cite{bn, jpn:nkg}, 
in the case where the non-linearity $G$ is given precisely by $G(x)=x|x|^2$. 
Although this is a special case of $(\ref{behaviour})$, the phenomenology 
of this particular equation is completely different from the general case,
as this non-linearity \emph{is} of the form $V'(x)$ for a potential $V\ge 0$. 
In view of the nonnegativity of $V$, the specific
form of $G$,\footnote{in particular, the sufficiently \emph{low} power} 
and the Sobolev embedding $H^1\hookrightarrow L^6$, 
global regularity for \emph{large data} solutions can be proven using
only the energy estimate arising from contraction of
the energy momentum tensor of $\phi$ with the Killing vector $\partial_t$.\footnote{These
methods have recently been extended~\cite{jpn:kerr} to the case of a Kerr background, where
$\partial_t$ is no longer everywhere timelike in the domain of outer communications.
In this case, global regularity is proven even though the solution is allowed to grow.}
In particular, no appeal need be made  
to decay for the linearized problem.
Under spherical symmetry, similar results can in fact be 
proven in the \emph{self-gravitating} case, i.e.~for the system
\begin{equation}
\label{sg1}
R_{\mu\nu}-\frac12g_{\mu\nu}R=2T_{\mu\nu}
\end{equation}
\begin{equation}
\label{sg2}
\Box_g\phi=V'(\phi)
\end{equation}
\begin{equation}
\label{sg3}
T_{\mu\nu}=\phi_{,\mu}\phi_{,\nu}-\frac12g_{\mu\nu}\phi^{,\lambda}\phi_{,\lambda}
-g_{\mu\nu}V(\phi),
\end{equation}
for a wide class of potentials.
Specifically, if the potential is bounded below,
i.e.~$V\ge -C>-\infty$, then a weak form of 
stability\footnote{i.e., that the solution does not blow up in the past of infinity}
has been proven in~\cite{md:higgs} for arbitrary
spherically symmetric solutions
of~\eqref{sg1}--\eqref{sg3}, in the case were it is assumed
\emph{a priori} that the past of infinity has a nontrivial complement. 
If $V\ge 0$, the results of~\cite{md:higgs} 
imply in particular that null infinity is
complete~\cite{md:sssts}.\footnote{In the spherically symmetric self-gravitating
case, an energy estimate is provided via the coupling by considering evolution equations
for the Hawking
mass $m$.} 
These questions received attention~\cite{hhm:gccv,hhm:update}
recently in the string theory community in the 
context of possible cosmic censorship violation.

In contrast to the above, for
the case of the ``John problem'' 
considered here, where we do \emph{not} assume a special structure for $G$,
results are only expected for small data solutions, and for $p$ sufficiently \emph{large}.
We recall the situation in Minkowski space:
The study of small amplitude solutions of the equation 
\begin{equation}\label{nlw-john}
\Box \phi = |\phi|^p
\end{equation}
on ${\mathbb R}^{3+1}$
was initiated by F.~John~\cite{john}, 
who showed that for $p>1+\sqrt 2$, solutions to~\eqref{nlw-john} 
 with sufficiently small, compactly supported initial data
remain regular for all time, and moreover, decay to $0$ 
along all inextendible causal geodesics. 
On the other 
hand, for $1<p<1+\sqrt 2$, F.~John established the 
existence of small data spherically symmetric solutions 
which blow up in finite time.\footnote{Consideration of the general $n$-dimensional problem 
leads to the critical power $p_c$ defined as the positive root of the quadratic equation 
$(n-1)p_c^2-(n+1)p_c-2=0$; this power first arose in work of W.~Strauss. 
For analogues of F.~John's result for ${\Bbb R}^{n+1}$,
see~\cite{gls} and references therein.} This blow-up result 
was later extended by J.~Schaeffer~\cite{schaef} to $p=1+\sqrt 2$. 

Recently, it was proven in~\cite{cg} that for
equation~\eqref{nlw-john} on a Schwarzschild background\footnote{The spherically symmetric
geometric wave operator on Schwarzschild can be reinterpreted, after multiplication of $\phi$
by $r$, as 
the ${\mathbb R}^{1+1}$ wave operator, with an additional positive potential,
and this is the approach taken in~\cite{cg}. For
extensions of the John problem to large classes of linear potentials, 
not necessarily positive, see~\cite{karageorgis}, and references therein.}
for $1< p<1+\sqrt{2}$, there exist solutions 
arising from arbitarily small compactly supported data
which blow up on the exterior region 
$J^-(\mathcal{I}^+)\cap J^+(\mathcal{I}^-)$.
In view of the restriction $p>4$ in our Theorem~\ref{main},
it would be interesting to understand the behaviour in the interval
$1+\sqrt{2}\le p \le 4$, in particular, to determine the 
critical power of $p$ separating stability and blow up.\footnote{For simplicity, we appeal 
in this paper only to the uniform $v$-decay
estimate--not the complete set of decay rates--proven
in \cite{mi:mazi}. We expect that a more detailed study
should lower the constraint $p>4$. It is not clear, however,
whether the estimates of~\cite{mi:mazi} alone will be able to retrieve the result for 
$p>1+\sqrt{2}$; for comparison with the case of ${\mathbb R}^{3+1}$,
see~\cite{john}.}

Besides addressing what is by now a classical problem in non-linear wave equations,
our motivation for studying this ``John problem'' in the black hole context
arises as follows:
Experience from non-linear stability problems
for the Einstein vacuum equations without symmetry strongly suggests 
that using decay rates of linearized fields will be an essential part of any argument. 
In the black hole context, implementing this approach 
introduces a number of technical challenges, as the notion of ``time'' with
respect to which decay is to be measured must necessarily be substantially altered. 
In this context, the present problem 
appears to concern the simplest geometrically well-defined
non-linear equation where ``global regularity'' is indeed shown via
\emph{decay} for the linearized problem. Moreover, we believe that 
the estimates proven here clarify the role of the event horizon--and 
the celebrated red-shift effect--for ensuring stability.  
Indeed, it is not clear whether the results shown here could be obtained
using only the (non-uniform) $t$-decay shown in~\cite{clsy:wpgs, st:pc}, which do not
precisely express behaviour up to the horizon. 
Further clarification of these issues may shed light on
the problem of the non-linear stability of the Kerr vacuum solution.

\section{Schwarzschild and Reissner-Nordstr\"om}
Let $(\mathcal{M},g)$ denote a 
non-extremal Reissner-Nordstr\"om solution corresponding
to parameters $M$, $e$, with $0\le|e|<M$.
(We refer the reader to standard references, for instance~\cite{he:lssst}.)
By convention here, the term Reissner-Nordstr\"om solution will always refer to a globally
hyperbolic spacetime which is the Cauchy development of
an asymptotically flat hypersurface $\Sigma$ with two ends\footnote{i.e., in
the $e\ne0$ case, we do not mean an analytic (or otherwise) extension thereof}.
Considering the quotient
\[
\mathcal{Q}=\mathcal{M}/SO(3),
\]
and defining the area-radius function
$r:\mathcal{Q}\to{\mathbb R}$
\begin{equation}
\label{arearadius}
r(q)=\sqrt{Area(q)/4\pi},
\end{equation}
we have that 
\[
g=\bar{g}+r^2 d\sigma_{{\mathbb S}^2},
\]
where $d\sigma_{{\Bbb S}^2}$ denotes the standard metric on the unit sphere, and
where $\bar{g}$ defines a $1+1$-dimensional Lorentzian metric on $\mathcal{Q}$.
The Lorentzian manifold $(\mathcal{Q}, \bar{g})$ can be covered
by a bounded global null coordinate system $(\bar{u},\bar{v})$. 
The range of such a coordinate system
is depicted, as a subset of $2$-dimensional Minkowski space,
in the Schwarzschild ($M>0$, $e=0$) case:
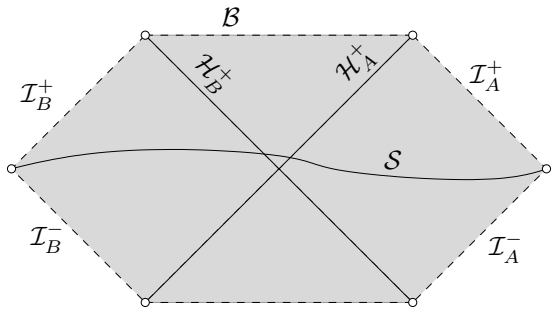
\[
\input{schwarzschild.pstex_t}
\]
and the charged Reissner-Nordstr\"om ($M>|e|>0$) case:
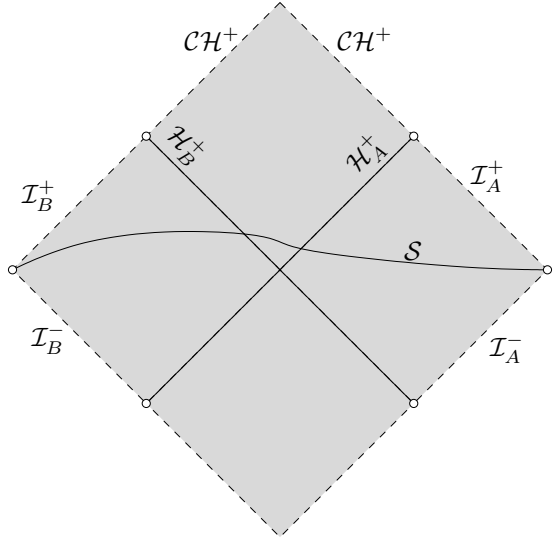
\[
\input{Re-No.pstex_t}
\]
These depictions are commonly known as \emph{Penrose diagrams}.
The curve $\mathcal{S}$ depicts the projection
to $\mathcal{Q}$ of a particular choice
of complete Cauchy surface $\Sigma\subset \mathcal{M}$.

We define the
\emph{Hawking mass} function $m:\mathcal{Q}\to{\mathbb R}$
by
\[
m=\frac r2\left(1-\bar{g}(\nabla r,\nabla r)\right),
\]
and the \emph{mass ratio} 
\[
\mu=\frac{2m}r.
\]
The mass $m$ is related to the parameter $M$ by
\[
M=m+\frac{e^2}{2r}.
\]
The sets $J^-(\mathcal{I}^+_A)\cap J^+(\mathcal{I}^-_A)$
and $J^-(\mathcal{I}^+_B)\cap J^+(\mathcal{I}^-_B)$, the two so-called
\emph{domains of outer communications}, can each be covered by a
coordinate system $(r,t)$ with coordinate range $(r_+,\infty)\times(-\infty,\infty)$, 
where $r$ is defined by $(\ref{arearadius})$,
and in which the metric $\bar{g}$ 
of $\mathcal{Q}$ takes the form
\[
\bar{g}=-\left(1-\frac{2M}r+\frac{e^2}{r^2}\right)dt^2+\left(1-\frac{2M}r+\frac{e^2}{r^2}\right)^{-1}
dr^2.
\]
Here $r_+$ is defined by
\begin{equation}
\label{r_+=}
r_+=M+\sqrt{M^2-e^2}.
\end{equation}
One immediately sees that the vector field $\partial_t$ is timelike Killing in these domains.

To define a null coordinate system, we first introduce the so-called
\emph{Regge-Wheeler} coordinate 
\begin{equation}
\label{Regge-Wheeler}
r_*=r+r_+^2(r_+-r_-)^{-1}\log|r-r_+|-r_-^2(r_+-r_-)^{-1}\log |r-r_-|,
\end{equation}
where
\[
r_-=M^2-\sqrt{M^2-e^2}.
\]
This coordinate satisfies
\[
\frac{d r_*}{d r}=(1-\mu)^{-1}
\]
We now set
\[
u=t-r_*,
\]
\[
v=t+r_*.
\]
With respect to these coordinates, the metric $\bar{g}$ takes the form
\[
\bar{g}=-(1-\mu)dudv=-\left(1-\frac{2M}r+\frac{e^2}{r^2}\right)dudv
\]
in $J^-({\mathcal I}^+_A)\cap J^+({\mathcal I}^-_A)$
and $J^-({\mathcal I}^+_B)\cap J^+({\mathcal I}^-_B)$; i.e.~$(u,v)$ are indeed
null. Individually, the coordinates $u$ and $v$ are each known as \emph{Eddington-Finkelstein}
retarded and advanced time coordinates, 
respectively.\footnote{The term \emph{Eddington-Finkelstein coordinate system}, however,
typically refers to the hybrids $(u,r)$ and $(v,r)$.}
Each of the domains of outer communications is
covered by the range $(-\infty,\infty)\times(-\infty,\infty)$.

The notion of asymptotic flatness is captured by the fact that
$1-\mu\to1$ as $r\to \infty$.
To examine the geometry of the horizon, first, we note that 
the timelike Killing vector field $\partial_t$ extends to a 
null tangent Killing vector field on the event 
horizon $\mathcal{H}=\mathcal{H}^+_A\cup\mathcal{H}^+_B$, vanishing 
at the sphere of bifurcation $\mathcal{H}^+_A\cap\mathcal{H}^+_{B}$.
Thus $1-\mu=0$ along $\mathcal{H}$, and moreover
$m$ and $r$ are constant. In terms of the
parameters, we have
\[
m=2r=2r_+=2(M+\sqrt{M^2-e^2}).
\]

This vanishing of $1-\mu$ 
relative to the Regge-Wheeler 
$(t,r_*)$ and null $(u,v)$ coordinates is captured by:
\begin{eqnarray}
\label{behaviour}
\nonumber
\hat{C}_1e^{\frac{(r_+-r_-)}{2r_+^2}(v-u)}&=&
\hat{C}_1e^{\frac{(r_+-r_-)}{r_+^2}r_*}\\
&\le&C_1(r-r_+)\\
\nonumber
&\le& 1-\mu\le C_2(r-r_+)\le  
\hat{C}_2e^{\frac{(r_+-r_-)}{r_+^2}r_*}=\hat{C}_2e^{\frac{(r_+-r_-)}{2r_+^2}(v-u)}
\end{eqnarray}
for small enough $-\infty<r^*$, and 
for some positive constants, $C_1$, $\hat{C}_1$, $C_2$, and $\hat{C}_2$, depending 
only on $M$ and $e$.

\section{Local existence}
Let $G$ be a bounded $C^\infty$ function of its argument, 
let $(\mathcal{M},g)$ be as in the previous section, and let
$\Sigma$ be a Cauchy surface.
We have the following
\begin{proposition}
\label{IVPstatement}
Let $\phi_0$ and $\phi_1$ be $C^\infty$ functions on
$\Sigma$. There exists
a unique globally hyperbolic subset $\mathcal{D}\subset\mathcal{M}$,
and a $C^\infty$ function $\phi:\mathcal{D}\to{\mathbb R}$ such that
\begin{enumerate}
\item
$\Sigma\subset \mathcal{D}$ is a Cauchy surface, 
and $\phi|_{\Sigma}=\phi_0$,
$N^\alpha\phi_\alpha=\phi_1$, where $N$ is the future-directed
unit normal to $\Sigma$.
\item
$\Box_g\phi=G(\phi)$
\item
The set $\mathcal{D}$ is maximal, i.e.~if $\mathcal{D}\subset\tilde{\mathcal{D}}$
is globally hyperbolic, with $\tilde\phi:\tilde{\mathcal{D}}\to {\mathbb R}$
satisfying the above, then $\tilde{\mathcal{D}}=\mathcal{D}$.
\end{enumerate}
We shall call $\mathcal{D}$ the \emph{maximal domain of existence}
of the solution.
Moreover, if $\Sigma$ is spherically symmetric, 
and $\phi_0$ and $\phi_1$ are constant 
on the $SO(3)$ orbits, then $\mathcal{D}$ is spherically
symmetric and $\phi$ is constant on the $SO(3)$ orbits,
and thus descends to a function on $\mathcal{D}_{\mathcal{Q}}
=\mathcal{D}/SO(3)$. 
\end{proposition}
In the spherically symmetric case, the regularity assumption
in the statement can be lowered to $G\in C^0_{loc}$, 
$\phi_0\in C^1_{loc}$, $\phi_1\in C^0_{loc}$.
The resulting solution $\phi$ is then $C^1$. In a slight abuse
of notation, we will denote $\mathcal{D}_{\mathcal{Q}}$ again as
$\mathcal{D}$. It turns out that
blowup can be easily characterized in terms of the $L^\infty$ norm 
of $\phi$,
i.e.~we have
\begin{proposition}
\label{ext}
Let $q\in\mathcal{Q}\cap\overline{\mathcal{D}}$. If 
$J^-(q)\cap \mathcal{D}\ne\emptyset$, and 
\[
\sup_{x\in J^-(q)\cap\mathcal{D}}|\phi(x)|<\infty,
\]
then $q\in \mathcal{D}$.
\end{proposition}
A similar statement can be made about the past.
The proofs of the above proposition are completely standard and are thus omitted.

\section{The main theorem}
The main theorem of this paper is the following:
\begin{theorem}
\label{main}
Let $\Sigma$ be a spherically symmetric Cauchy surface
in a Reissner-Nordstr\"om spacetime $(\mathcal{M},g)$, with parameters 
$0\le e<M$, and let $\phi_0$, $\phi_1$ be $C^1$, $C^0$ functions, respectively, on $\Sigma$,
constant
on the $SO(3)$ orbits, and supported in $\Sigma\cap\{r\le R\}$, for some $R<\infty$. 
Let $\mathcal{D}$ denote the maximal domain of existence for
the initial value problem defined by Proposition~\ref{IVPstatement},
where $G$ is assumed to satisfy $(\ref{satisfies})$ with $p>4$. 
Then there exists an $\epsilon$, depending
only on $M$, $e$, $\Sigma$, and $R$, such that
for 
\[
|\phi_0|_{C^1}<\epsilon,\footnote{The norm above  refers to the induced Riemannian
metric on $\Sigma$.}
\]
\[
|\phi_1|_{C^0}<\epsilon,
\]
$\mathcal{D}$ contains the closures of the two 
domains of outer communications,
\[
\overline{(J^-(\mathcal{I}^+_A)\cap 
J^+(\mathcal{I}^-_A))\cup (J^-(\mathcal{I}^+_B)\cap J^+(\mathcal{I}^-_B))}
\subset\mathcal{D}.
\] 
Moreover,
\[
|\phi(u,v)|\le C (\max\{v,1\})^{-1},
\]
and, for each fixed $r_0>r_+$,
\[
|\phi(r,t)|\le C(r_0)(|t|+1)^{-1}\qquad r\ge r_0,
\]
where $v$ and $t$ denote an Eddington-Finkelstein advanced time coordinate
and a static time coordinate, respectively,
on either of the two domains of outer communications.
\end{theorem}

\section{Proof of Theorem~\ref{main}: The bootstrap}
In view of simple arguments of Cauchy stability, and the fact that the
data are compactly supported on $\Sigma$, it follows, that given any
$t_0$, say $t=1$, then
for $\epsilon$ sufficiently small, we have
\[
\{t=1\}\subset \mathcal{D}
\]
and
\begin{equation}
\label{newsmallness1}
|\phi|_{t=1}|_{C^1}\le \tilde{\epsilon},
\end{equation}
\begin{equation}
\label{newsmallness2}
|N\phi|_{t=1}|_{C^0}\le\tilde{\epsilon}.
\end{equation}
Here again, the norms refer to the induced Riemannian metric of
$\{t=1\}$, and $N$ refers to a unit future directed normal.
The relation of $\mathcal{S}=\Sigma/SO(3)$ and $\{t=1\} $ could be
given, for instance, as follows:
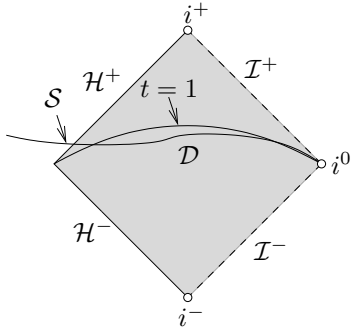
\[
\input{ext.pstex_t}
\]
Note $\phi|_{t=1}$ will \emph{not} in general be compactly supported
in $(r_+,\infty)$. It will be supported in $(r_+, \bar R]$ for some $\bar R<\infty$.
Moreover, we have $\tilde{\epsilon}=\tilde{\epsilon}(\epsilon,\Sigma, \bar R, M, e)\to0$, 
as $\epsilon\to0$.


In proving our theorem, it clearly suffices to show
that 
\[
\mathcal{D}\supset J^+(\Sigma)\cap \overline{J^-(\mathcal{I}^+_A)\cap J^+(\mathcal{I}^-_A)}.
\]
Define the
set $\mathcal{B}\subset \mathcal{D}\cap J^-(\mathcal{I}^+_A)\cap J^+(\mathcal{I}^-_A)$
by
\begin{equation}
\label{Tdef}
\mathcal{B}=\left\{q\in \mathcal{D}\cap J^-(\mathcal{I}^+_A)\cap J^+(\mathcal{I}^-_A)
:\sup_{x\in J^-(q)\cap\{t\ge 1\}\cap
\mathcal{D}}v_+
(x)|\phi(x)|\le B\right\}
\end{equation}
where $v_+$ denotes $\max\{v,1\}$,
and for a constant $B$ to be determined later.

This subset is clearly non-empty and open for 
$B>\tilde\epsilon v_+(R,1)$.
In view of Proposition~\ref{ext}, proving Theorem~\ref{main} reduces to 
showing that $\mathcal{B}$ is closed. This would follow immediately from 
\begin{theorem}
\label{equiv}
For $q\in \mathcal{B}$, we have 
\[
\sup_{x\in J^-(q)\cap\{t\ge 1\}\cap
\mathcal{D}}v_+(x)\left|\phi(x)\right|\le \frac B2.
\]
\end{theorem}

To prove this, 
we shall view our non-linear problem in $\mathcal{B}$
as a solution of 
\begin{equation}\label{eq:lin-sch}
\Box_g \psi =F,
\end{equation}
where $F$ is to satisfy the estimates arising from $(\ref{Tdef})$ if
$F$ is set equal to the non-linearity.
We thus turn in the next section to the study of this linear inhomogeneous equation.
The proof of Theorem~\ref{equiv} will then follow easily in Section~\ref{theproofof}.

\section{The linear problem}
First some formulas:
Written with respect to Regge-Wheeler coordinates, 
equation 
\eqref{eq:lin-sch} takes the form
\begin{equation}\label{eq:lin-reg}
(1-\mu)^{-1} \Big (\pa_t^2\psi - r^{-2}\pa_{r_*} (r^2 \pa_{r_*} \psi)\Big ) = F.
\end{equation}
Alternatively,  relative to the null coordinates $(u=t-r_*, v=t+r_*)$
equation \eqref{eq:lin-sch} can be written as a first order system
\begin{align}
&\pa_v (r\pa_u \psi) = \frac  {(1-\mu) }{2r} r\pa_v\psi + \frac{(1-\mu)}4 r F,
\label{eq:lin-u}\\
&\pa_u (r\pa_v \psi) = \frac  {(1-\mu) }{2r} r\pa_u\psi + \frac{(1-\mu)}4 r F
\label{eq:lin-v}
\end{align}
for the pair of functions $\theta=r\pa_v \psi$ and $\zeta=r\pa_u \psi$, or
by a single equation
\begin{equation}\label{eq:lin-psi}
\pa_v \pa_u  (r\psi) = -\left(M-\frac{e^2}r\right)
\frac  {(1-\mu) }{2r^2} \psi + \frac{(1-\mu)}4 r F.
\end{equation}
We also record the following form of equation \eqref{eq:lin-u}, crucial 
to understanding the so-called \emph{red-shift effect}
\begin{equation}\label{eq:lin-red}
\pa_v \frac{ r \pa_u \psi}{1-\mu} + \frac {1}{r^2}
 \left(M-\frac{2e^2}r\right)
 \frac{r\pa_u \psi}{1-\mu}
=  \frac  {r\pa_v\psi}{2r} + \frac 14 r F.
\end{equation}
This effect has been discussed in~\cite{mi:mazi}. 

From \eqref{eq:lin-reg} and \eqref{eq:lin-u}--\eqref{eq:lin-v}
we can derive the following energy estimates.
Given $(u_1,v_1)$, for $t\le (v_1+u_1)/2$, 
define
\begin{eqnarray*}
E_t^{(u_1,v_1)}[\psi]&=&
\int_{t-u_1}^{-t+v_1}\big (|\pa_t \psi(t,r_*)|^2 +
|\pa_{r_*}\psi(t,r_*)|^2\big ) r^2 \,dr_*\\
&=&2 \int_{t-u_1}^{-t+v_1} \big (|r\pa_u \psi(t,r_*)|^2 +|r\pa_{v}\psi(t,r_*)|^2\big )  \,dr_*.
\end{eqnarray*}
For $t_1=(v_1+u_1)/2\ge t_2\ge t_3$, we have
the following energy identity:
\begin{align*}
E_{t_2}^{(u_1,v_1)}[\psi]&+\int_{2t_3-u_1}^{2t_2-u_1} |r\pa_v\psi(u_1,v)|^2\,dv+
\int_{2t_3-v_1}^{2t_2-v_1} |r\pa_u\psi(u,v_1)|^2\,du\\ &=
E_{t_3}^{(u_1,v_1)}[\psi]
+
\frac 12\int_{t_3}^{t_2} \int_{t-u_1}^{-t+v_1}
(1-\mu) F(t,r_*) \,\pa_t\psi(t,r_*)\, r^2\,dr_*\,dt.
\end{align*}
The estimate refers thus to the region depicted below:
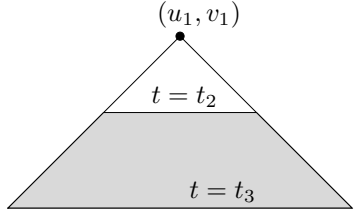
\[
\input{eest.pstex_t}
\]
We also have the following characteristic energy
identity in the characteristic rectangle $[u_2,u_1]\times [v_2,v_1]$:
\begin{align*}
\int_{u_2}^{u_1}|r\pa_u\psi(u,v_1)|^2  du &+\int_{v_2}^{v_1}
|r\pa_{v}\psi(u_1,v)|^2  \,dv=\\
&\int_{u_2}^{u_1}|r\pa_u\psi(u,v_0)|^2   du +\int_{v_2}^{v_1}
|r\pa_{v}\psi(u_2,v)|^2 \,dv\\ 
&+\frac 14\int_{u_2}^{u_1} \int_{v_2}^{v_1}
(1-\mu) F(u,v) (\pa_u\psi(u,v) - \pa_v \psi(u,v))\, r^2\,dv\,du.
\end{align*}
The reader can refer to the diagram of Theorem~\ref{pricelaw1}.

\subsection{The homogeneous problem}
The linear homogeneous problem, i.e.~the case $F=0$, 
was studied in~\cite{mi:mazi},
where pointwise decay estimates were proven. 
In this section, we will present
these estimates in a form suitable for understanding the inhomogeneous problem
via Duhamel's principle.


First we have the following preliminary result:
\begin{theorem} 
\label{pricelaw0}
Fix $r_0>r_+$, and let $(u_1,v_1)$ be an arbitrary point. 
Let $t_0<(u_1+v_1)/2$, and consider the lighter-shaded region
\[
\mathcal{E}=J^-(r_1,t_1)\cap J^+(\{t=t_0\})\cap\{r\ge r_0\}
\]
depicted below\footnote{Note that despite the diagram,
we are not necessarily assuming that $r(u_1,v_1)\ge r_0$.}:
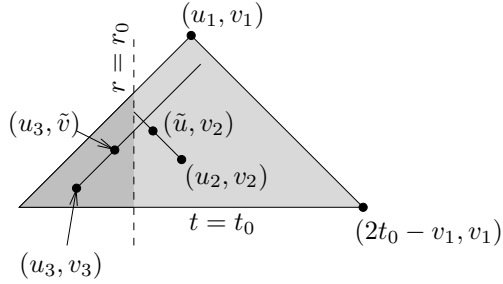
\[
\input{prelim.pstex_t}
\]
Let $\psi$ be a solution to
\[
\Box_g\psi=0
\]
on $J^-(u_1,v_1)\cap J^+(\{t=t_0\})$.
Let $(u_2,v_2)\in\mathcal{E}$.
Then there exists a constant $A$, depending only on $r_0$, $e$, and $M$,
such that 
\begin{equation}
\label{easiest}
r|\psi(\tilde{u},v_2)|^2\le 
A\left (r|\psi(u_2,v_2)|^2+E^{(\tilde{u},v_2)}
_{t_0}[\psi]\right)  
\end{equation}
\begin{equation}
\label{easiest2}
r^3|\partial_v\psi(\tilde{u},v_2)|^2\le A
\left(r^3|\partial_v\psi(u_2,v_2)|^2
+ E^{(\tilde{u},v_2)}_{t_0}[\psi]\right),
\end{equation}
for any point $\tilde{u}\ge u_2$ with $(\tilde{u},v_2)\in \mathcal{E}$.

On the other hand, if $(u_3,v_3)$ is any point in 
$J^-(u_1,v_1)\cap J^+(\{t=t_0\})$,
then
\begin{equation}\label{eq:mu}
\left(\frac{r|\partial_u\psi(u_3,\tilde{v})|}{1-\mu}\right)^2\le A\left( 
e^{-\frac{r_+-r_-}{r^2(u_3,\tilde v)}(\tilde v-v_3)}
\left(\frac{r|\partial_u\psi(u_3,v_3)|}{1-\mu}\right)^2
+E^{(u_3,\tilde{v})}_{t_0}[\psi]\right)
 \end{equation}
for any $\tilde{v}\ge v_3$ 
such that $(u_3,\tilde{v})\in  J^-(r_1,t_1)\cap J^+(\{t=t_0\})$.

Finally, if 
\[
C=\sup_{\mathcal{E}\cap \{r\ge r_0\}}
\left\{|r\psi|(r,t_0), |r^2\partial_v(r\psi)(r,t_0)|\right\}
\]
then
\begin{equation}\label{eq:nu}
|r\psi(r,t)|^2+ |r^2\partial_v\psi(r,t)|^2+|r^2\partial_v(r\psi)(r,t)|^2\le 
A\left(C^2+E^{(u_1,v_1)}_{t_0}[\psi]\right),
\end{equation}
throughout $J^-(r_1,t_1)\cap\{r\ge r_0\}\cap \{t\ge t_0\}$.
\end{theorem}
\begin{proof}
Equations~\eqref{easiest}, \eqref{easiest2}, and \eqref{eq:nu} 
can be easily deduced from the equations 
\eqref{eq:lin-u}--\eqref{eq:lin-psi} and the energy 
estimates of the previous section. For the estimate \eqref{eq:mu},
we write the homogeneous wave equation in the red-shift
from  given by \eqref{eq:lin-red}:
\[
\pa_v  \frac{r\pa_u\psi}{1-\mu} + \frac 1{r^2}
\left(M-\frac {e^2}r\right)  \frac{r\pa_u\psi}{1-\mu}=
\frac 1{2r}{r\pa_v\psi}.
\]
Introducing the integrating factor $e^{\int_{v_3}^v\frac 1{r^2}(M-\frac{e^2}r)dv}$
and integrating in $v$ we obtain
$$
\frac{r\pa_u\psi(u_3,\tilde v)}{1-\mu}=
e^{-\int_{v_3}^{\tilde v} \frac 1{r^2}(M-\frac{e^2}r)dv} 
\,\,\frac{r\pa_u\psi(u_3,v_3)}{1-\mu}+
\int_{v_3}^{\tilde v}e^{-\int_{v}^{\tilde v} \frac 1{r^2}(M-\frac{e^2}r)dv'}
\,\,\frac 1{2r}{r\pa_v\psi}\, dv.
$$
Since 
\begin{equation}
\label{surfgrav}
M-\frac {e^2}r\ge M-\frac {e^2}{r_+}=\frac{r_+-r_-}2>0,
\end{equation}
 and $r$ is increasing in the $v$-direction, we conclude
\begin{eqnarray*}
\frac{r|\pa_u\psi(u_3,\tilde v)|}{1-\mu}&\le& 
e^{ -\frac{r_+-r_-}{2r^2(u_3,\tilde v)}(\tilde v-v_3)} 
\,\,\frac{r|\pa_u\psi(u_3,v_3)|}{1-\mu}\\
&&\hbox{}+
\int_{v_3}^{\tilde v} e^{ -\frac{r_+-r_-}{2r^2(u_3,\tilde v)}(v-v_3)} 
\,\,\frac 1{2r}{r|\pa_v\psi(u_3,v)|}\, dv .
\end{eqnarray*}
Applying Cauchy-Schwarz gives
$$
\frac{r|\pa_u\psi(u_3,\tilde v)|}{1-\mu}\le 
e^{ -\frac{r_+-r_-}{2r^2(u_3,\tilde v)}(\tilde v-v_3)} 
\,\,\frac{r|\pa_u\psi(u_3,v_3)|}{1-\mu}+ A
\big (\int_{v_3}^{\tilde v} \,{|r\pa_v\psi(u_3,v)|^2}\, dv \big)^{\frac 12}.
$$
\end{proof}
We now turn to decay estimates. The full results of~\cite{mi:mazi} are 
complicated to state, as different regions must be treated separately. 
In this paper we shall only use the following:
\begin{theorem}
\label{pricelaw1}
Let $(u_1, v_1)$, $(u_2,v_2)$ be points 
with $(u_1,v_1)\in J^+(u_2,v_2)$, i.e.~with $u_1\ge u_2$, $v_1\ge v_2$.
Consider the characteristic rectangle 
\[
\mathcal{R}= J^+(u_2,v_2)\cap J^-(u_1,v_1)
\]
depicted below:
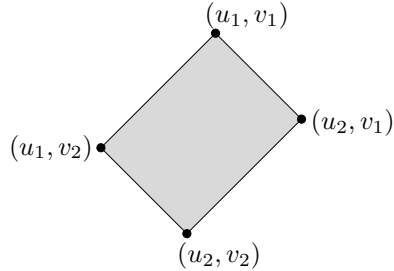
\[
\input{price.pstex_t}
\]
Let $\psi$ be a solution to
$\Box_g\psi=0$ on $\mathcal{R}$.
Let
\begin{eqnarray*}
C=\max\left\{
\right.&\sup_{v_2\le v\le v_1}
|r\partial_v\psi(u_2,v)|+|r^2\partial_v(r\psi)(u_2,v)|,\\
&\hbox{\ }\left.
\sup_{u_2\le u\le u_1}
|(1-\mu)^{-1}\partial_u\psi(u,v_2)|\right\}.
\end{eqnarray*}
Then there exists a constant $A$ depending only on 
$r_0=r(u_2,v_2)$, $e$ and $M$
such that 
\[
\psi(u_1,v_1)\le AC(v_1-v_2)^{-1}.
\]
\end{theorem}

The above theorem concerns a characteristic initial value problem.
We can combine Theorems \ref{pricelaw0} and \ref{pricelaw1} to derive 
decay estimates for solutions of the homogeneous Cauchy problem.
\begin{theorem}\label{pricelaw2}
Let $(u_1,v_1)$ be an arbitrary point, and
and let $t_0<(v_1+u_1)/2$. 
Let $\psi$ be a solution to
$$
\Box_g\psi=0 
$$
on $J^-(u_1,v_1)\cap J^+(\{t=t_0\})$.
Let
\begin{equation}
\label{C=}
C=\sup_{t=t_0}
\{(1-\mu)^{-\frac 12}r^2|\partial_t\psi|,(1-\mu)^{\frac 12}r^2|\partial_r\psi|,
r^2|\partial_v (r\psi)|, r|\psi|\}.
\end{equation}
Pick a point $(u_2,v_2)\in J^-(u_1,v_1)\cap\{t=t_0\}$:
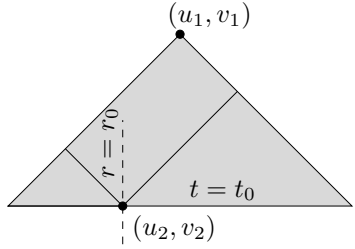
\[
\input{homog.pstex_t}
\]
It follows that
\begin{equation}
\label{yIldIz}
|\psi(u_1,v_1)|\le AC(v_1-v_2)^{-1},
\end{equation}
where $A$ depends only on $r_0=r(u_2,v_2)$, $e$ and $M$.
\end{theorem}
\begin{proof}
We first show that $E^{(u_1,v_1)}_{t_0}[\psi]\le A C^2$. Note that
$\pa_{r_*}\psi=(1-\mu)\pa_r\psi$.
\begin{align*}
\left(E^{(u_1,v_1)}_{t_0}[\psi]\right)^2 
=&\int_{t_0-u_1}^{-t_0+v_1} r^2 (|\partial_t\psi|^2+|\partial_{r_*}\psi|^2)\, dr_*\\
=&\int_{ r(r_*=t_0-u_1)}^{r(r_*=-t_0+v_1)} 
(1-\mu)^{-1} r^2 (|\partial_t\psi|^2+|\partial_{r_*}\psi|^2)\,dr
\\ \le& \ C^2\int_{r_+}^{\infty} r^{-2}\, dr 
 \le C^2 r_+^{-1}.
\end{align*}

We apply now Theorem~\ref{pricelaw0}.
From~\eqref{eq:mu} and~\eqref{C=}, the estimate 
\begin{eqnarray}
\label{puttogether}
\nonumber
|(1-\mu)^{-1} \pa_u\psi(u,v_2)|&\le& e^{-\frac{(r_+-r_-)}{r^2(u,v_2)}(v_2-v_0(u))}
|(1-\mu)^{-1} \pa_u\psi(u,v_0(u))|+ AC\\
				&\le& e^{-\frac{(r_+-r_-)}{r^2(u,v_2)}(v_2-v_0(u))}
(1-\mu)^{-\frac12}C+AC
\end{eqnarray}
follows for 
$u_2\le u\le u_1$, where $v_0(u)$ denotes the $v$-coordinate of the 
intersection of the outgoing constant-$u$ ray with $t=t_0$. 
Since,
$$
v_2-v_0=2(r_{0*}-r_*(u,v_2)),
$$
we have by~\eqref{behaviour}, for sufficiently small $r_*=r_*(u,v_2)$, $r=r(u,v_2)$:
\begin{eqnarray}
\label{puttogether2}
\nonumber
(1-\mu(r_*))^{-\frac12}e^{-\frac{r_+-r_-}{2r^2}(v_2-v_0)}
	&\le& 
			A(1-\mu(r_*))^{-\frac12}e^{\frac{r_+-r_-}{r^2}r_*}\\
	\nonumber
	&\le& A	 e^{-\frac{r_+-r_-}{2r_+^2}r_*+\frac{r_+-r_-}{r^2}r_*}\\
	\nonumber
	& = & A	 e^{-\frac{(r-\sqrt{2}r_+)(r+\sqrt{2}r_+)(r_+-r_-)}{2r_+^2r^2}r_*}\\
	\nonumber
	&\le& A e^{-r_*(r-\sqrt{2}r_+)}\\
	&\le& A.
\end{eqnarray}
(Note that by convention here, $A$ does not always represent the same constant.)
Since the above estimate is also clearly true for large $r$, putting together
\eqref{puttogether} and \eqref{puttogether2} yields
\begin{equation}\label{eq:mu-1}
\sup_{u_2\le u\le u_1} (1-\mu)^{-1} |\pa_u\psi(u,v_2)|\le AC.
\end{equation}

On the other hand, from~\eqref{eq:nu},
we have
\[
|\pa_v\psi(u_2,v)|\le ACr^{-1},
\] 
\[
|\partial_v(r\psi)(u_2,v)|\le ACr^{-2},
\]
for $v_1\le v\le v_2$, in view
of the fact that $r\ge r_0$ for such $(u_0,v)$. Our result now
follows immediately from Theorem~\ref{pricelaw1}.
\end{proof}
\begin{remark}
Observe that the constant $C$ in the above theorem is dominated by 
the (Riemannian)
$C^1$ norm of $r^3\psi$,
$$
C\le |r^3\psi|_{t=t_0}|_{C^1}.
$$ 
\end{remark}
\begin{remark}
The quantity ${\pa_{\bar{u}} \psi}/{\pa_{\bar{u}} r}$
is clearly independent of the choice of retarded time $\bar{u}$ and
coincides with $-(1-\mu)^{-1} \pa_u \psi$
for Eddington-Finkelstein $u$. As one can choose a 
global regular null coordinate system,
such that $\pa_{\bar{u}} r\le -1$ along the ray $u=u_2$ in a 
neighborhood of $\mathcal{H}^+_A\cap\{u=u_2\}$, 
one could have alternatively derived 
\eqref{eq:mu-1} by appealing to Cauchy stability.
\end{remark}
 

\subsection{The inhomogeneous problem}
To prove Theorem~\ref{equiv},
we must derive estimates for solutions of the equation
\begin{equation}
\label{savgrammikn}
\Box_g\Psi=F,
\end{equation}
with vanishing initial data at $t=1$.
Let
$\psi(r,t;s)$
denote the solution of
$\Box_g\psi=0$
evaluated at $(r,t)$
with
initial condition
\begin{equation}
\label{ic}
\psi(\cdot,s)=0,\qquad 
\pa_t \psi(\cdot,s)=(1-\mu)F(\cdot,s).
\end{equation}
Using representation \eqref{eq:lin-sch} we can 
then write the solution of $(\ref{savgrammikn})$ as
\begin{equation}
\label{integral}
\Psi(r,t)=\int_1^t\psi(r,t;s)ds.
\end{equation}
We prove the following theorem:
\begin{theorem}\label{pricelaw4}
Let $(u_1,v_1)$ be given, 
let $\Psi$ satisfy $(\ref{savgrammikn})$ on $J^-(u_1,v_1)\cap
\{t\ge 1\}$, with $\Psi=0$, $\partial_t\Psi=0$ on $J^-(u_1,v_1)\cap\{t=1\}$,
and let $\alpha>4$.
Then there exists a constant $A$, depending
only on $M$, $e$, and $\alpha$, such that
\begin{equation}\label{eq:Psi-total-}
|\Psi(u_1,v_1)| \le  A  (v_1)^{-1} 
\sup_{J^-(r_1,t_1)\cap\{t\ge 1\}}
 v_+^\alpha|F(u,v)|,
\end{equation}
where $v_+=\max (v,1)$.
\end{theorem}
\begin{proof} 
Let $r_1=r(u_1,v_1)$, $t_1=(v_1+u_1)/2$.
We will fix some $r_0>r_+$, to be determined later,
and distinguish the cases $r_1\ge r_0$ and $r_1<r_0$.

\subsubsection{Estimates in the region $r\ge r_0$}
Consider first the case $r_1\ge r_0$.
We have $v_1=t_1+r_1^*\ge t_1+r_0^*$.

Let us take $s> 2v_1^\omega$, for $0\le \omega<1$.
To compute $\psi(r_1,t_1;s)$, we decompose 
\begin{align*}
J^-(r_1,t_1)\cap \{t=s\} &= \{(u,v): \,\, u=2s-v,\,\,v_1-2(t_1-s)\le v\le v_1\}\\ 
&=I_1\cup I_2
\end{align*}
as depicted:
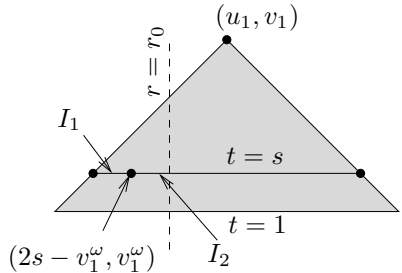
\[
\input{far.pstex_t}
\]
i.e.~we have
\begin{align*}
I_1&= \{(u,v): \,\, u=2s-v,\,\,v_1-2(t_1-s)\le v\le v_1^\omega\},\\
I_2&= \{(u,v): \,\, u=2s-v,\,\,v_1^\omega\le v\le v_1\}.
\end{align*}
Note that on the set $I_1$,
$$
2r_*= v-u = 2(v-s)  \le  - 2 v_1^\omega  
$$
and therefore, by $(\ref{behaviour})$, we have
$$
1-\mu(r,s)\le \hat{C}_2 e^{-\frac{(r_+-r_-)}{r_+^2}v_1^\omega},\qquad
(r,s)\in I_1.
$$
By the above estimate, \eqref{easiest} and \eqref{ic}, we have 
\begin{align}
\label{above-est}
\nonumber
|\psi(r_1,t_1;s) |^2&\le
 A^2\left(E_s[\psi(\cdot,\cdot;s)]\right)^2\\
 \nonumber
 &= A^2 \int_{I_1\cup I_2}
 (1-\mu) |F(r_*,s)|^2 r^2\, dr\\ &\le 
 \tilde{A}^2\left(e^{-\frac{(r_+-r_-)}{r_+^2}v_1^\omega} \sup_{I_1} |F(r_*,s)|^2 + 
 v_1^3 \sup_{I_2} |F(r_*,s)|^2\right),
\end{align}
where in the last inequality we used that $r\le C'v_1$ on $I_2$.

\vskip 2pc

On the other hand, for $s\le 2v_1^\omega$ we can apply 
Theorem \ref{pricelaw2} for the points 
$(r_1,t_1)=(u_1,v_1)$ and $(r_0,s)=(u_2,v_2)$, where 
\[
v_2=s+ r_0^*\le 2v_1^\omega.
\]
(In the last inequality, we assumed that $r_0$ is chosen in particular so that
$r_{0*}\le 0$.)  Then,
\begin{align*}
|\psi(r_1,t_1;s)|&\le A v_1^{-1} \sup_r (1-\mu)^{\frac 12} r^3 |F(r,s)|
\\ &\le A v_1^{-1} \sup_{r^*\le -\frac s2} (1-\mu)^{\frac 12} r^3 |F(r,s)| +
A v_1^{-1} \sup_{r^*\ge -\frac s2} (1-\mu)^{\frac 12} r^3 |F(r,s)|.
\end{align*}
For points $r_*\le -\frac s2$, estimate \eqref{behaviour} gives 
$$
1-\mu\le \hat{C}_2e^{-\frac{(r_+-r_-)}{2r_+^2} s},
$$
while for $r_*\ge -\frac s2$,
we clearly have
$$
v=r_*+s\ge \frac s2.
$$
Thus, for $\delta>0$,
\begin{eqnarray}
\label{taalla}
\nonumber
|\psi(r_1,t_1;s)|
&\le& A v_1^{-1} e^{-\frac{(r_+-r_-)}{2r_+^2} s} \sup_{r_*\le -\frac s2} r^3 |F(r,s)|\\
&&\hbox{} +
A v_1^{-1} s^{-1-\delta}\sup_{r_*\ge -\frac s2} v^{1+\delta} r^3 |F(r,s)|.
\end{eqnarray}

Applying $(\ref{integral})$, and the estimates $(\ref{above-est})$ and $(\ref{taalla})$,
we obtain
\begin{align*}
|\Psi(r_1,t_1)|\le& t_1A e^{-\frac{(r_+-r_-)}{2r_+^2}v_1^\omega} \sup_{J^-(r_1,t_1)} 
|F(r,s)| +A  t_1v_1^{\frac 32} \sup_{J^-(r_1,t_1)\cap \{v\ge v_1^\omega\}}
|F(u,v)| \\ &+ A v_1^{-1} \sup_{J^-(r_1,t_1)}  v_+^{1+\delta} r^3 |F(r,s)|\\
\le& (v_1)_+^{\frac 52} \sup_{J^-(r_1,t_1)}
|F(u,v)| +(v_1)_+^{-1} \sup_{J^-(r_1,t_1)}  v_+^{1+\delta} r^3 |F(r,s)|,
\end{align*}
where we have used the fact that for $r_1\ge r_0>r_+$, 
we have $(v_1)_+ \ge\sim  t_1$, with constant depending on $r_0$.

Clearly, there exists a constant $H$, depending only on $M$ and $e$, such that for
$r\ge H$, 
we have $v=t+r_*\ge r $ in the region $t\ge 1$. It follows immediately that for all $r$, we have
$r\le H +|v|$ in the region $t\ge 1$. We obtain, for any $\alpha\ge 0$,
\begin{align}
|\Psi(r_1,t_1)|\le& 
A  (v_1)_+^{\frac 52-\omega\alpha} \sup_{J^-(r_1,t_1)\cap \{v\ge v_1^\omega\}}
 v_+^\alpha |F(u,v)|\nonumber \\ 
 &+
A(v_1)_+^{-1} \sup_{J^-(r_1,t_1)}  v_+^{4+\delta} |F(u,v)|\\
\le&
A  (v_1)_+^{\frac 52-\omega\alpha} \sup_{J^-(r_1,t_1)}
 v_+^\alpha |F(u,v)|\nonumber \\ 
 &+
A (v_1)_+^{-1} \sup_{J^-(r_1,t_1)}  v_+^{4+\delta} |F(u,v)|.\nonumber
\end{align}
Choosing 
$$
\alpha=4+\delta,\qquad \omega=\frac 7{2(4+\delta)}
$$
we obtain
\begin{equation}\label{eq:Psi-big}
|\Psi(r_1,t_1)|\le A (v_1)_+^{-1} \sup_{J^-(r_1,t_1)}  v_+^\alpha |F(u,v)|.
\end{equation}
\subsubsection{Estimates in the region $r<r_0$.}
\label{redshiftsec}
Consider now the case $r_+<r_1< r_0$. In this region we shall use 
the 
\emph{red-shift effect} 
manifest\footnote{See~\cite{mi:mazi}. The point, as we shall see,
is that we can extract decay in $v$ for $\partial_u\Psi$, even though
we are integrating in $v$ towards the future.} 
in the equation \eqref{eq:lin-red}:
\[
\pa_v \frac{ r \pa_u \Psi}{1-\mu} +\frac 1{r^2}\left(M-\frac{e^2}r\right)
\frac{r\pa_u \Psi}{1-\mu}
=  \frac  {r\pa_v\Psi}{2r} + \frac 14 r F.
\]

For fixed $u$, let $v_*(u)$ be defined by $t(u,v_*(u))=1$.
Using the integrating factor $e^{\int_{v_*(u)}^v r^{-2}(M-{e^2}r^{-1})\, dv}$ 
and the vanishing 
of the initial data on $\{t=1\}$, we obtain
\[
\frac{ r \pa_u \Psi (u,v)}{1-\mu} = 
\frac 14 \int_{v_*(u)}^v
e^{-\int_{v'}^v r^{-2}(M-{e^2}r^{-1})} \Big (2\pa_v \Psi(u,v') + r F(u,v')\Big )\, dv'
\]
Integrating by parts we derive 
\begin{align*}
&\frac{ r \pa_u \Psi (u,v)}{1-\mu} =  \frac 12 \Psi (u,v)\\
&+
\int_{v_*(u)}^v
e^{-\int_{v'}^v r^{-2}(M-{e^2}r^{-1})} \Big (-2r^{-2}(M-e^2r^{-1})
\Psi(u,v') + r F(u,v')\Big )\, dv'.
\end{align*}
Since in the region under consideration, we have the bound
$r_+< r\le r_0$, and also, 
\[
M-\frac{e^2}r>M-\frac{e^2}{M+\sqrt{M^2-e^2}}=\sqrt{M^2-e^2}>0,
\]
it follows 
that
$$
\frac{ |\pa_u \Psi (u,v)|}{1-\mu} \le   A|\Psi (u,v)| + A\int_{v_*(u)}^v
e^{-c(v-v')} \Big (|\Psi(u,v')| + | F(u,v')|\Big )\, dv'.
$$
Given $v$, let $u_*(v)$ be defined by $r(u_*(v),v)=r_0$. Then 
\begin{align*}
|\Psi (u,v)|&\le |\Psi(u_*(v),v)|  + A\int_{u_*(v)}^u 
(1-\mu(u',v)) |\Psi(u',v)|\, du' \\ &+  A\int_{u_*(v)}^u\int_{v_*(u')}^v
e^{-c(v-v')} \Big (|\Psi(u',v')| + | F(u',v')|\Big )\, dv' (1-\mu(u',v))\,du'.
\end{align*}
We now observe that in the region $r_+<r\le r_0$,
we can estimate 
$|1-\mu |< \delta$ for some $\delta$ that can be made
arbitrarily small provided that $r_0$ is sufficiently 
close to $r_+$.  Moreover, 
\[
\int_{u_*(v)}^u |1-\mu(u',v)|\, du' = r(u,v) - r_0 <\delta.
\]
Therefore, for any $u\ge u_*(v)$,
\begin{eqnarray}\label{eq:Psi}
 |\Psi(u,v)|&\le& |\Psi(u_*(v),v)| + 
 A\delta\sup_{u'\in [u_*(v),u]}{\Psi(u',v)}\\
  \nonumber
 &&\hbox{}+
 A\delta 
 \sup_{u'\in [u_*(v),u]}\int_{v_*(u')}^v
e^{-c(v-v')} \Big (|\Psi(u',v')| + | F(u',v')|\Big )\, dv'
\end{eqnarray}
Split 
\begin{align*}
\int_{v_*(u')}^v
e^{-c(v-v')} &\Big (|\psi(u',v')| + | F(u',v')|\Big )\, dv'= \\
&\int_{v_*(u')}^{v/2}
e^{-c(v-v')} \Big (|\psi(u',v')| + | F(u',v')|\Big )\, dv'\\ &+
\int_{v/2}^v
e^{-c(v-v')} \Big (|\psi(u',v')| + | F(u',v')|\Big )\, dv'
\end{align*}
and estimate 
\begin{align*}
\int_{v_*(u')}^{v/2}
e^{-c(v-v')}& \Big (|\Psi(u',v')| + | F(u',v')|\Big )\, dv'\\
&\le
A e^{-c v/2} \sup_{v'\in [v_*(u'),v]  }\Big (|\Psi(u',v')| + | F(u',v')|\Big ),\\
\int_{v/2}^v
e^{-c(v-v')}& \Big (|\Psi(u',v')| + | F(u',v')|\Big )\, dv'\\
&\le
A \sup_{v'\in [v/2,v] }\Big (|\Psi(u',v')| + | F(u',v')|\Big ).
\end{align*}
Returning to \eqref{eq:Psi}, we derive 
\begin{align*}
 |\Psi(u,v)| \le& 
 |\Psi(u_*(v),v)| +
 A\delta\sup_{u'\in [u_*(v),u]}{\Psi(u',v)}\\
&+
 A\delta \sup_{u'\in [u_*(v),u]}
\sup_{v'\in [v/2,v]  }\Big (|\Psi(u',v')|+ | F(u',v')|\Big )\\
 &+A\delta  e^{-c v} 
\sup_{u'\in [u_*(v),u]}
\sup_{v'\in [v_*(u'),v]  }\Big (|\Psi(u',v')| + | F(u',v')|\Big ).
\end{align*}
Thus,
\begin{align*}
|v_+^\beta\Psi(u,v)|
\le& 
 v_+^\beta|\Psi(u_*(v),v)| +
 A\delta v_+^\beta\sup_{u'\in [u_*(v),u]}{\Psi(u',v)}\\
 &
 +v_+^\beta A\delta \sup_{u'\in [u_*(v),u]}
\sup_{v'\in [v/2,v]  }|\Psi(u',v')|\\
&+Av_+^\beta \delta  e^{-c v} 
\sup_{u'\in [u_*(v),u]}
\sup_{v'\in [v_*(u'),v]  }|\Psi(u',v')|\\
&+v_+^\beta A\delta \sup_{u'\in [u_*(v),u]}
\sup_{v'\in [v/2,v]  }| F(u',v')|\\
&+Av_+^\beta \delta  e^{-c v} 
\sup_{u'\in [u_*(v),u]}
\sup_{v'\in [v_*(u'),v]  } | F(u',v')|\\
\le& 
 v_+^\beta|\Psi(u_*(v),v)| 
+ A\delta v_+^\beta\sup_{u'\in [u_*(v),u]}{\Psi(u',v)}\\
 &+ A\delta \sup_{u'\in [u_*(v),u]}
\sup_{v'\in [v/2,v]  }|(v'_+)^\beta\Psi(u',v')|\\
&+A \delta
\sup_{u'\in [u_*(v),u]}
\sup_{v'\in [v_*(u'),v]  }|\Psi(u',v')|\\
&+ A\delta \sup_{u'\in [u_*(v),u]}
\sup_{v'\in [v/2,v]  }(v'_+)^\beta| F(u',v')|\\
 &+A \delta  
\sup_{u'\in [u_*(v),u]}
\sup_{v'\in [v_*(u'),v]  }| F(u',v')|\\
\le&
 v_+^\beta|\Psi(u_*(v),v)| 
  +A\delta v_+^\beta\sup_{u'\in [u_*(v),u]}{\Psi(u',v)}\\
&+  2A\delta \sup_{u'\in [u_*(v),u]}
\sup_{v'\in [v_*(u'),v]  }|(v'_+)^\beta\Psi(u',v')|\\
&+
2A\delta \sup_{u'\in [u_*(v),u]}
\sup_{v'\in [v_*(u'),v]  }(v'_+)^\beta| F(u',v')|.
\end{align*}
The constant $A$ is independent of the choice of $r_0$ and thus
of $\delta$.
\emph{A fortiori}, the above estimate applies
also when $|v_+^\beta\Psi(u,v)|$ on the left hand side is replaced by
$\sup_{u'\in [u_*(v),u]}
\sup_{v'\in [v_*(u'),v]  }|(v'_+)^\beta\Psi(u',v')|$. Thus, if
$r_0$ is chosen close enough to $r_+$, so that $\delta$ is sufficiently
small, we have at our point $(r_1,t_1)$ the estimate
\begin{align}
 |\Psi(r_1,t_1)| \le 2
(v_1)_+^{-\beta} \Big (& \sup_{v\in [v^*,v_1]  } v_+^\beta |\Psi(u_*(v),v)| 
\nonumber\\ &+ \sup_{u\in [u_*(v_1),u_1]}
\sup_{v\in [v_*(u),v_1]}   v_+^\beta | F(u,v)|\Big ),\label{eq:Psi-les}
\end{align}
for any $\beta\ge 0$.
Here $v^*$ denotes the $v$-coordinate of the intersection of the curve 
$r=r_0$ and $\{t=1\}$. Observe that the first term on the right hand side of
the above inequality can be rewritten as 
$\sup_{1\le t\le t_1 } v_+^\beta|\Psi(r_0,t)|$;
the supremum is thus taken over $\Psi$ evaluated at points
in a subset of the
region
$r\ge r_0$. We choose 
\[
\beta=1.
\] 

From \eqref{eq:Psi-big} and \eqref{eq:Psi-les}, we infer now that
\begin{equation}\label{eq:Psi-total}
|\Psi(r_1,t_1)| \le  A  (v_1)_+^{-1} 
 \max_{J^-(r_1,t_1)}  v_+^\alpha |F(u,v)|,
\end{equation}
for \emph{all} $(r_1,t_1)$.
\end{proof}

\section{Proof of Theorem~\ref{equiv}}
\label{theproofof}
Recall the set $\mathcal{B}$ and the constant $B$ defined by \eqref{Tdef}.
Let $(r_1,t_1)=(u_1,v_1)\in \mathcal{B}$.
Set $F=G(\phi)$ and consider 
$\Psi$ a solution of $(\ref{savgrammikn})$ 
in $J^-(u_1,v_1)$ with vanishing data at $t=1$,
and let $\psi$ be a solution of the homogeneous equation with
data
\[
\psi|_{J^-(u_1,v_1)\cap \{t=1\}}=\phi|_{J^-(u_1,v_1)\cap \{t=1\}},
\]
\[
\partial_t\psi|_{J^-(u_1,v_1)\cap \{t=1\}}=\partial_t\phi|_{J^-(u_1,v_1)\cap \{t=1\}}.
\]
We have that 
\begin{equation}
\label{a9roisma}
\phi=\psi+\Psi
\end{equation}
in $J^-(u_1,v_1)\cap J^+(\{t=1\})$.

By Theorems~\ref{pricelaw0} and~\ref{pricelaw1}, 
it follows that
\begin{equation}
\label{homogeneouspsi}
|\psi(u_1,v_1)| \le E (v_1)_+^{-1},
\end{equation}
where $E=E(\epsilon,M,e,\Sigma, R)\to0$ as $\epsilon\to0$.
If $0\le B<c$, then
by assumption~\eqref{satisfies},
\[
|G(\phi(u,v))|\le K|\phi(u,v)|^p, \forall (u,v)\in \mathcal{B}
\]
Thus, by the definition $(\ref{Tdef})$ 
of the set $\mathcal{B}$, we have
\[
|F(u,v)|\le KB^pv_+^{-p}.
\]
Therefore for $p> 4$, choosing $\alpha=p$ and applying Theorem~\ref{pricelaw4},
we have
\begin{equation}
\label{proproteleutaio}
|\Psi(u_1,v_1)|\le A(v_1)_+^{-1}KB^p
\end{equation}
and thus, by $(\ref{a9roisma})$ and $(\ref{homogeneouspsi})$, 
\begin{equation}
\label{proteleutaio}
|\psi(u_1,v_1)|\le A(v_1)_+^{-1}KB^p+E(v_1)_+^{-1}.
\end{equation}
Let us define $B= 4E$,
and let us require $E$ (and thus $\epsilon$) to 
be sufficiently small
so that $AKB^p< \frac B4$.
Inequality $(\ref{proteleutaio})$ then yields 
\[
|\psi(u_1,v_1)|\le \frac B2(v_1)_+^{-1},
\]
as required.

\section{Acknowledgements}
M.~D.~is supported in part by NSF grant DMS-0302748. I.~R.~is
supported in part by NSF grant DMS-0406627. Part of this research
was conducted when I.~R.~was a long-term prize fellow
of the Clay Mathematics Institute.

\end{document}